\documentclass[12pt, twoside, epsf]{article}

\usepackage{color,graphicx,times}
\usepackage{amsmath, amssymb, graphics}

\newcommand{\mathsym}[1]{{}}
\newcommand{\unicode}[1]{{}}

\makeatletter
\long\def\M#1{\leavevmode\setbox\@tempboxa\hbox{#1}\@tempdima\fboxrule
    \advance\@tempdima \fboxsep \advance\@tempdima \dp\@tempboxa
   \hbox{\lower \@tempdima\hbox
  {\vbox{\hrule \@height \fboxrule
          \hbox{  \hskip\fboxsep
          \vbox{\vskip\fboxsep \box\@tempboxa\vskip\fboxsep}\hskip
                 \fboxsep\vrule \@width \fboxrule}%
                  }}}}
\makeatother

\let \ttorg \tt \def \tt{\ttorg \obeyspaces}

\begin{document}

\date{}

\title{\bf Majorana Fermions and Representations of the Braid Group}

\author{Louis H. Kauffman \\
  Department of Mathematics, Statistics and Computer Science \\
  University of Illinois at Chicago \\
  851 South Morgan Street\\
  Chicago, IL, 60607-7045}

\maketitle
  
\thispagestyle{empty}

\begin{abstract}In this paper we study unitary braid group representations associated with
Majorana Fermions. Majorana Fermions are represented by {\it Majorana operators}, elements of a Clifford algebra. The paper recalls and proves a general result about braid group representations associated with Clifford algebras, and
compares this result with the Ivanov braiding associated with Majorana operators.  The paper generalizes observations of Kauffman and Lomonaco and of Mo-Lin Ge to show that certain strings of Majorana operators give
rise to extraspecial 2-groups and to braiding representations of the Ivanov type. 
\end{abstract}

\noindent {\bf Keywords.~~}{knots, links, braids, braid group, Fermion, Majorana Fermion, Kitaev chain, extraspecial 2-group, Majorana string, Yang-Baxter equation, quantum process, quantum computing.}\\

\noindent {\bf Mathematics Subject Classification}{ 2010: 57M25, 81Q80.}

\section{Introduction}
 In this paper we study a Clifford algebra generated by non-commuting elements of square equal to one and the relationship of this algebra with braid group representations and Majorana Fermions \cite{Ivanov,Kitaev,KnotLogic}.  
 Majorana Fermions can be seen not only in the structure of collectivities of electrons, as in the quantum Hall effect \cite{Wilczek}, but also in the structure of single electrons both by experiments with electrons in nanowires \cite{Kouwenhouven,Beenakker} and also by the decomposition of the operator algebra for a Fermion \cite{Hatfield}  into a Clifford algebra generated by two Majorana operators. 
Majorana Fermions have been discussed by this author and his collaborators in \cite{Kitaev,KnotLogic,KLogic,AnyonicTop,SpinTop,KL,BG,Fibonacci,BMF,Iterants,IAlg,RK} and parts of the present paper depend strongly on this previous work. In order to make this paper 
 self-contained, we have deliberately taken explanations, definitions and formulations from these previous papers, indicating the references when it is appropriate.
The purpose of this paper is to discuss these  braiding representations, important for relationships among physics, quantum information and topology. \\

We study Clifford algebra and its relationship with braid group representations related to Majorana Fermion operators. Majorana Fermion operators $a$ and $b$  are defined, so that the creation and annihilation operators  $\psi^{\dagger}$ and $\psi$ for a single standard Fermion can be expressed through them. The Fermion operators satisfy the well known algebraic rules: $$(\psi^{\dagger})^{2} = \psi^{2} = 0,$$ 
 $$\psi \psi^{\dagger} + \psi^{\dagger} \psi = 1.$$ Remarkably, these equations are satisfied if we take $$\psi = (a + i b)/2,$$ $$\psi^{\dagger} = (a - ib)/2$$ where the Majorana operators $a,b$ satisfy $$a^{\dagger} = a, b^{\dagger} = b,$$ $$a^2 = b^2 = 1, ab + ba = 0.$$ In certain situations, it has been conjectured and partially verified by experiments \cite{Kouwenhouven,Kitaev} that electrons (in low temperature nano-wires) may behave as though each electron were physically a pair of Majorana particles described by these Majorana operators. In this case the 
mathematics of the braid group representations that we study may have physical reality.\\

Particles corresponding to the Clifford algebra generated by $a$ and $b$ described in the last paragraph are called Majorana particles because they satisfy $a^{\dagger} = a$ and $b^{\dagger} = b$ indicating that they are 
their own anti-particles. Majorana \cite{Majorana} analyzed real solutions to the Dirac equation \cite{D,KN1} and conjectured the existence of such particles that would be their own anti-particle. It has been conjectured that the neutrino is such a particle. Only more recently \cite{Ivanov,Kitaev} has it been suggested that electrons may be composed of pairs of Majorana particles. It is common to speak of Majorana particles when referring to particles that
satisfy the interaction rules for the original Majorana particles. These interaction rules are, for a given particle $P$, that $P$ can interact with another identical $P$ to produce a single $P$ or to produce an annihilation. For this, we write $PP = P + 1$ where the right hand side is to be read as a superposition of the possibilities $P$ and $1$ where $1$ stands for the state of annihilation, the absence of the particle $P.$ We refer to this equation as
{\it the fusion rules for a Majorana Fermion}. In modeling the quantum Hall effect \cite{Wilczek,Fradkin,B1,B2}, the braiding of  quasi-particles (collective excitations) leads to non-trival
representations of the Artin braid group. Such particles are called {\it anyons}. The braiding in these models is related to 
topological quantum field theory. \\

Thus there are two algebraic descriptions for Majorana Fermions -- the fusion rules,  and the associated Clifford algebra. One may use both the Clifford algebra and the fusion rules in a single physical situation. However, for studying braiding, it turns out that the Clifford algebra leads to braiding and so does the fusion algebra in the so-called Fibonacci model (while the Fibonacci model is not directly related to the Clifford algebra). Thus we can discuss these two forms of braiding. We show mathematical commonality between them in the appendix to the present paper. Both forms of braiding could be present in a single physical system. For example,
in the quantum Hall systems, the anyons (collective excitiations of electrons)  can behave according to the Fibonacci model, and the edge effects of these anyons can be modeled using Clifford algebraic braiding for Majorana Fermions.\\

Braiding operators associated with Majorana operators are described as follows. Let $\{c_1 , c_2 , \cdots , c_n \}$ denote a collection of Majorana operators such that $c_{k}^2 = 1$ for $k=1, \cdots, n$ and
$c_{i}c_{j} + c_{j}c_{i} = 0$ when $i \ne j.$ Take the indices $\{ 1,2,...,n\}$ as a set of residuces modulo $n$ so that $n+1 = 1.$ Define operators $$\sigma_{k} = (1 + c_{k+1}c_{k})/\sqrt{2}$$ for $k= 1,\cdots n$ where it is understood that $c_{n+1} = c_{1}$ since $n +1 = 1$ modulo $n.$ Then one can verify that 
$$\sigma_{i}\sigma_{j} = \sigma_{j}\sigma_{i}$$ when $|i - j| \ge 2$ and that $$\sigma_{i}\sigma_{i+1}\sigma_{i} = \sigma_{i+1}\sigma_{i}\sigma_{i+1}$$ for all $i=1,\cdots n.$ Thus $$\{\sigma_1, \cdots , \sigma_{n-1} \}$$
describes a representation of the $n$-strand Artin braid group $B_{n}.$ As we shall see in Section 3, this representation has very interesting properties and it leads to unitary representations of the braid group that can support partial topological computing. What is missing to support full topological quantum computing in this representation is a sufficient structure of $U(2)$ transformations. These must be supplied along with the braiding 
operators. It remains to be seen if the braiding of Majorana operator constituents of electrons can be measured, and if the physical world will yield this form of partial topological computing.\\

Here is an outline of the contents of the paper.
Section 2 reviews the definition of the Artin braid group and emphasizes that the braid group on $n$ strands, $B_{n},$ is a natural extension of the symmetric group on $n$ letters. Along with its topological interpretations, this proximity to the symmetric group probably explains the many appearances of the braid group in physical and mathematical problems. Section 3 dicusses how unitary braiding operators can be (in the presence of local unitary transformations) universal gates for quantum computing and why certain solutions of the Yang-Baxter (braiding) equation, when entangling, are such universal gates. Section 4 discusses how a Clifford algebra of 
Majorana Fermion operators can produce braiding representations. We recall the Clifford Braiding Theorem of \cite{BMF} and we show how extraspecial 2-groups give rise to representations of the Artin Braid Group.
In Section 4 we show how Majorana Fermions can be used to construct represenations of the Temperley-Lieb algebra and we analyze corresponding braid group representations, showing that they are equivalent to our
already-constructed representations from the Clifford Braiding Theorem. In section 6 we show how the Ivanov \cite{Ivanov} representation of the braid group on a space of Majorana Fermions generates a $4 \times 4$ universal quantum
gate that is also a braiding operator. This shows how Majorana Fermions can appear at the base of (partial) topological quantum computing. We say partial here because a universal topological gate of this type must be supported by local unitary transformations not necessarily generated by the Majorana Fermions. In section 7 we consider the  Bell Basis Change Matrix $B_{II}$  and braid group representations that are related to it.
The matrix itself is a solution to the Yang-Baxter equation and so it is a universal gate for partial topological computing \cite{BG}. Mo-Lin Ge has observed that $B_{II} = (I + M)/\sqrt{2}$ where $M^2 = -I.$ In fact, we take
\begin{equation}
M=\left[\begin{array}{cccc}
0 & 0 & 0 & 1\\
0 & 0 & -1 & 0\\
0 & 1 & 0 & 0\\
-1 & 0 & 0 & 0
\end{array}\right] =
\left[\begin{array}{cccc}
1 & 0 & 0 & 0\\
0 & -1 & 0 & 0\\
0 & 0 & 1 & 0\\
0 & 0 & 0 & -1
\end{array}\right]
\left[\begin{array}{cccc}
0 & 0 & 0 & 1\\
0 & 0 & 1 & 0\\
0 & 1 & 0 & 0\\
1 & 0 & 0 & 0
\end{array}\right].
\end{equation}
Let
\begin{equation}
A =
\left[\begin{array}{cccc}
1 & 0 & 0 & 0\\
0 & -1 & 0 & 0\\
0 & 0 & 1 & 0\\
0 & 0 & 0 & -1
\end{array}\right], 
B= \left[\begin{array}{cccc}
0 & 0 & 0 & 1\\
0 & 0 & 1 & 0\\
0 & 1 & 0 & 0\\
1 & 0 & 0 & 0
\end{array}\right].
\end{equation}
Thus $M= AB$ and $A^2 = B^2 = I$ while $AB + BA = 0.$ Thus we can take $A$ and $B$ themselves as Majorana Fermion operators. This is our key observation in this paper. The fact that the matrix $M$
factors into a product of (Clifford algebraic) Majorana Fermion operators means that the extraspecial 2-group associated with $M$ can be seen as the result of a {\it string} of Majorana operators where we 
define such a string as follows: A list of pairs $A_{k}, B_{k}$  of Majorana operators that satisfies the identities below is said to be a {\it string of Majorana operators.}

\begin{eqnarray}
A_{i}^2 =B_{i}^2 = 1,\\
A_{i}B_{i} = - B_{i}A_{i},\\
A_{i}B_{i+1} = - B_{i+1}A_{i},\\
A_{i+1}B_{i} = B_{i}A_{i+1},\\
A_{i}B_{j} = B_{j}A_{i}, \,\, for \,\, |i-j|>1 \\
A_{i}A_{j} = A_{j}A_{i} \\
B_{i}B_{j} = B_{j}B_{i} \,\, for \, all \,\, i \,\, and \,\, j.
\end{eqnarray}

Letting $M_{i} = A_{i}B_{i}$ we obtain an extraspecial 2-group and by the Extraspecial 2-Group Braiding Theorem of Section 4, a representation of the Artin Braid group with generators
$\sigma_{i} = (I + A_{i}B_{i})/\sqrt{2}.$ If we think of the Majorana Fermions $A_{i}$ and $B_{i}$ as anyonic particles, then it will be of great interest to formulate a Hamiltonian for them and to analyze them in analogy with 
the Kitaev spin chain. This project will be carried out in a sequel to the present paper.\\

In this paper, we continue a description of the relationship of the Majorana string for the $B_{II}$ matrix and the work of Mo-Lin Ge \cite{Ge} relating with an extraspecial 2-group and with the topological order in the 
Kitaev spin chain. We believe that our formulation of the Majorana string sheds new light on this relationship.\\

The last section of the paper is an appendix on braid group representations of the quaternions. As the reader can see, the different braid group representations studied in this paper all stem from the same formal considerations that occur in classifying quaternionic braid representations. The essential reason for this is that given three Majorana Fermions $A,B,C$ that pairwise anti-commute and such that $A^2 = B^2 = C^2 = 1,$ then
the Clifford elements $BA, CB, AC$ generate a copy of the quaternions. Thus our Clifford braiding representations are generalizations of particular quaternionic braiding. The appendix is designed to show that the quaternions are a rich source of connection among these representations, including the Fibonacci model \cite{AnyonicTop,Fibonacci}, which is briefly discussed herein.\\ 

\noindent {\bf Acknowledgement.}  
Much of this paper is based upon our joint work in the papers \cite{TEQE,Spie,Teleport,QK1,QK2,QK3,QK4,QK5,BG,AnyonicTop,QCJP1,QCJP2,Fibonacci,KLogic}. We have woven  this work into the present paper in a form that is coupled with recent and previous work on relations with logic and with Majorana Fermions.  It gives the author great pleasure to thank the Simons Foundation for support under Collaboration Grant, Award Number 426075.\\

\section{Braids}

\begin{figure}
     \begin{center}
     \begin{tabular}{c}
     \includegraphics[height=6cm]{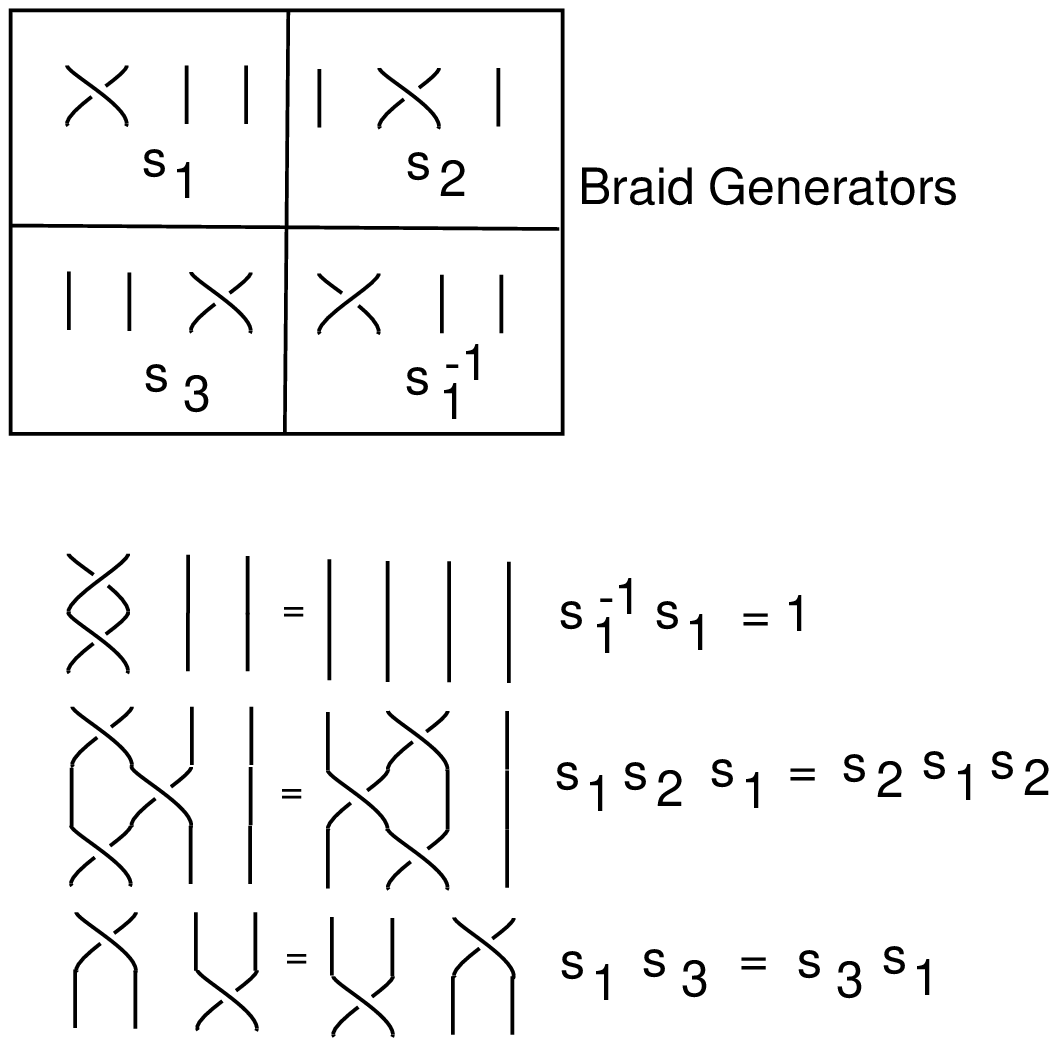}
     \end{tabular}
     \end{center}
     \caption{\bf Braid Generators  }
     \label{Figure 4 }
     \end{figure}

\begin{figure}
     \begin{center}
     \begin{tabular}{c}
     \includegraphics[height=6cm]{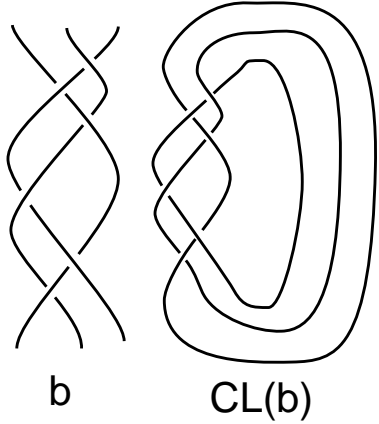}
     \end{tabular}
     \end{center}
     \caption{\bf Borromean Rings as a Braid Closure  }
     \label{Figure 6 }
     \end{figure}

A {\it braid} is an embedding of a collection of strands that have 
their ends in two rows of points that are set one above the other with respect to a choice of vertical. The strands are not
individually knotted and they are disjoint from one another. See Figures ~\ref{Figure 4 }, and ~\ref{Figure 6 } for illustrations of braids and moves on braids. Braids can be 
multiplied by attaching the bottom row of one braid to the top row of the other braid. Taken up to ambient isotopy, fixing the endpoints, the braids form
a group under this notion of multiplication. In Figure~\ref{Figure 4 } we illustrate the form of the basic generators of the braid group, and the form of the
relations among these generators. Figure~\ref{Figure 6 } illustrates how to close a braid by attaching the top strands to the bottom strands by a collection of 
parallel arcs. A key theorem of Alexander states that every knot or link can be represented as a closed braid. Thus the theory of braids is critical to the 
theory of knots and links. Figure~\ref{Figure 6 } illustrates the famous Borromean Rings (a link of three unknotted loops such that any two of the loops are unlinked)
as the closure of a braid.\\

Let $B_{n}$ denote the Artin braid group on $n$ strands.
We recall here that $B_{n}$ is generated by elementary braids $\{ s_{1}, \ldots ,s_{n-1} \}$
with relations 

\begin{enumerate}
\item $s_{i} s_{j} = s_{j} s_{i}$ for $|i-j| > 1$, 
\item $s_{i} s_{i+1} s_{i} = s_{i+1} s_{i} s_{i+1}$ for $i= 1, \ldots, n-2.$
\end{enumerate}

     \bigbreak

\noindent See Figure~\ref{Figure 4 } for an illustration of the elementary braids and their relations. Note that the braid group has a diagrammatic
topological interpretation, where a braid is an intertwining of strands that lead from one set of $n$ points to another set of $n$ points.
The braid generators $s_i$ are represented by diagrams where the $i$-th and $(i + 1)$-th strands wind around one another by a single 
half-twist (the sense of this turn is shown in Figure~\ref{Figure 4 }) and all other strands drop straight to the bottom. Braids are diagrammed
vertically as in Figure~\ref{Figure 4 }, and the products are taken in order from top to bottom. The product of two braid diagrams is accomplished by
adjoining the top strands of one braid to the bottom strands of the other braid. \\

In Figure~\ref{Figure 4 } we have restricted the illustration to the
four-stranded braid group $B_4.$ In that figure the three braid generators of $B_4$ are shown, and then the inverse of the
first generator is drawn. Following this, one sees the identities $s_{1} s_{1}^{-1} = 1$ 
(where the identity element in $B_{4}$ consists in  four vertical strands), 
$s_{1} s_{2} s_{1} = s_{2} s_{1}s_{2},$ and finally
$s_1 s_3 = s_3 s_1.$ \\

Braids are a key structure in mathematics. It is not just that they are a collection of groups with a vivid topological interpretation.
From the algebraic point of view the braid groups $B_{n}$ are important extensions of the symmetric groups $S_{n}.$ Recall that the 
symmetric group $S_{n}$ of all permutations of $n$ distinct objects has presentation as shown below.\\

\begin{enumerate}
\item $s_{i}^{2} = 1$ for  $i= 1, \ldots n-1,$
\item $s_{i} s_{j} = s_{j} s_{i}$ for $|i-j| > 1$, 
\item $s_{i} s_{i+1} s_{i} = s_{i+1} s_{i} s_{i+1}$ for $i= 1, \ldots n-2.$
\end{enumerate}

Thus $S_{n}$ is obtained from $B_{n}$ by setting the square of each braiding generator equal to one. We have an exact sequence of groups
$${1} \longrightarrow P(n) \longrightarrow B_{n} \longrightarrow S_{n} \longrightarrow {1}$$ exhibiting the Artin braid group as an extension of the symmetric group.
The kernel $P(n)$ is the {\it pure braid group}, consisting in those braids where each strand returns to its original position. \\

In the next sections we shall show how representations of the Artin braid group, rich enough to provide a dense set of transformations in the 
unitary groups (see \cite{AnyonicTop} and references therein), arise in relation to Fermions and Majorana Fermions. Braid groups are {\it in principle} fundamental to quantum computation and quantum information theory.
\bigbreak

     \bigbreak

\section{Braiding Operators and Universal Quantum Gates}
 
A key concept in the construction of quantum link invariants is
the association of a Yang-Baxter operator $R$ to each elementary crossing in a
link diagram. The operator $R$ is a linear mapping  
$$R\colon \ V\otimes V \longrightarrow V\otimes V$$ 
defined on the  $2$-fold tensor product of a vector space $V,$ generalizing the permutation of the factors
(i.e., generalizing a swap gate when $V$ represents one qubit). Such transformations are not 
necessarily unitary in 
topological applications. It is useful to understand
when they can be replaced by unitary transformations 
for the purpose of quantum 
computing. Such unitary $R$-matrices can be used to 
make unitary representations of the Artin braid group.
\bigbreak

More information about the material sketched in this section can be found in \cite{KP,AnyonicTop,BG}.\\
 
\begin{figure}
     \begin{center}
     \begin{tabular}{c}
     \includegraphics[height=3cm]{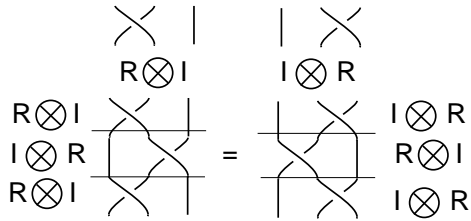}
     \end{tabular}
     \end{center}
     \caption{\bf The Yang-Baxter equation }
     \label{Figure 7 }
     \end{figure}

A solution to the Yang-Baxter equation, as described in the last 
paragraph is a matrix $R,$ regarded as a mapping of a
two-fold tensor product of a vector space
$V \otimes V$ to itself that satisfies the equation 

$$(R \otimes I)(I \otimes R)(R \otimes I) = 
(I \otimes R)(R \otimes I)(I \otimes R).$$ From the point of view of topology, the matrix $R$ 
is regarded as representing an elementary bit of braiding 
represented by one string
crossing over another. In Figure~\ref{Figure 7 } we have illustrated 
the braiding identity that corresponds to the Yang-Baxter equation.
Each braiding picture with its three input lines (below) 
and output lines (above) corresponds to a mapping of the three fold
tensor product of the vector space $V$ to itself, as required 
by the algebraic equation quoted above. The pattern of placement of the 
crossings in the diagram corresponds to the factors 
$R \otimes I$ and $I \otimes R.$ This crucial 
topological move has an algebraic
expression in terms of such a matrix $R.$  We need to study solutions of the Yang-Baxter equation that are unitary.
Then the $R$ matrix can be seen {\em either} as a braiding matrix 
{\em or} as a quantum gate in a quantum computer.

     \bigbreak

 \subsection{Universal Gates}
A {\em two-qubit gate} $G$ is a unitary  linear mapping $G:V \otimes V \longrightarrow V$ where $V$ is a two complex dimensional
vector space. We say that the gate $G$ is {\em universal for quantum computation} (or just {\em universal}) if $G$ together with 
local unitary transformations (unitary transformations from $V$ to $V$) generates all unitary transformations of the complex vector
space of dimension $2^{n}$ to itself. It is well-known \cite{N} that $CNOT$ is a universal gate. (On the standard basis,
$CNOT$ is the identity when the first qubit is $|0 \rangle $, and it flips the second qbit, leaving the first alone, when the first qubit is $|1 \rangle .$)
\bigbreak

\noindent A gate $G$, as above, is said to be {\em entangling} if there is a vector  
$$| \alpha \beta \rangle = | \alpha \rangle \otimes | \beta \rangle \in V \otimes V$$ such that 
$G | \alpha \beta \rangle$ is not decomposable as a tensor product of two qubits. Under these circumstances, one says that 
$G | \alpha \beta \rangle$ is {\em entangled}.
\bigbreak

\noindent In \cite{BB},  the Brylinskis
give a general criterion of $G$ to be universal. They prove that {\em a two-qubit gate $G$ is universal if and only if it is
entangling.} 
\bigbreak

\noindent {\bf Remark.} A two-qubit pure state $$|\phi \rangle = a|00 \rangle + b|01 \rangle + c|10 \rangle + d|11 \rangle$$
is entangled exactly when $(ad-bc) \ne 0.$ It is easy to use this fact to check when a specific matrix is, or is not, entangling.
\bigbreak

\noindent {\bf Remark.} There are many gates other than $CNOT$ that can be used as universal gates in the presence of local unitary 
transformations (see \cite{BG}). Some of these are themselves topological (unitary solutions to the Yang-Baxter equation, see \cite{BG,BA}) and themselves generate
representations of the Artin braid group. Replacing $CNOT$ by a solution to the Yang-Baxter equation does not place the local unitary transformations as
part of the corresponding representation of the braid group. Thus such substitutions give only a partial solution to creating topological 
quantum computation. In a full solution (e.g. as in \cite{AnyonicTop}) all unitary operations can be built directly from the braid group representations, and one hopes that the topology in the physics behind these
representations will give the system protection against decoherence. It remains to be seen if partial topological computing systems can have this sort of protection.\\

\section {\bf Fermions, Majorana Fermions and Braiding}
\noindent{\bf Fermion Algebra.} Recall Fermion algebra \cite{Hatfield}. One has Fermion annihilation operators $\psi$ and their
conjugate creation operators $\psi^{\dagger}.$
One has $\psi^{2} = 0 = (\psi^{\dagger})^{2.}$
There is a fundamental commutation relation
$$\psi \psi^{\dagger} + \psi^{\dagger} \psi = 1.$$
If you have more than one standard Fermion operator, say $\psi$ and $\phi$,
then they anti-commute:
$$\psi \phi = - \phi \psi.$$
The Majorana Fermions \cite{Majorana,Ivanov} $c$ satisfy $c^{\dagger} = c$ so that they
are their own anti-particles. They have a different algebra structure than standard Fermions, as we shall see below.
The reader may be curious just why the Majorana Fermions are assigned a Clifford algebra structure, as we are about to do.
We give a mathematical motivation below, by showing that, with this Clifford algebra structure, a standard Fermion is generated by two Majorana Fermions.
One way of putting this is that we can use $c_1=( \psi + \psi^{\dagger})/2$,
$c_2 = (\psi - \psi^{\dagger})/(2i)$ to make two operators satisfying $c_{i}^{\dagger} = c_{i}$ for $i = 1,2.$ We shall see that $c_{1}$ and $c_{2}$ satisfy Clifford algebra identities.\\

Majorana Fermion operators can model 
quasi-particles, and they are related to braiding and to topological
quantum computing.  A group of researchers  \cite{Kouwenhouven,Beenakker}
have found quasiparticle Majorana Fermions in edge effects in nano-wires.
(A line of Fermions could have a Majorana Fermion happen non-locally from
one end of the line to the other.) The Fibonacci model that we discuss \cite{AnyonicTop,Kitaev} is also based on
Majorana particles, possibly related to collective electronic excitations. 
If $P$ is a Majorana Fermion particle, then $P$ can interact with itself to either produce itself or to annihilate itself. This is the simple ``fusion algebra" for this particle. One can write
$P^2 = P + *$ to denote the two possible self-interactions of the particle $P,$ as we have discussed in the introduction.
The patterns of interaction and braiding of such a particle $P$ give
rise to the Fibonacci model \cite{AnyonicTop}.\\

\noindent {\bf Majoranas make Fermions.} Majorana operators \cite{Majorana} are related to standard Fermions as follows:
We first take two Majorana operators $c_{1}$ and $c_{2}$.
The algebra for Majoranas is $c_{1} = c_{1}^{\dagger}$, $c_{2} = c_{2}^{\dagger}$ and $c_{1}c_{2} = -c_{2}c_{1}$ if $c_{1}$ and $c_{2}$ are
distinct Majorana Fermions with  $c_{1}^{2}= 1$ and  $c_{2}^{2}= 1.$ Thus the operator algebra for a collection of Majorana particles is a Clifford algebra.
One can make a standard Fermion operator from two Majorana operators via
$$\psi = (c_{1} + ic_{2})/2,$$
$$\psi^{\dagger} = (c_{1} -ic_{2})/2.$$
Note, for example, that 
$$\psi^2 =  (c_{1} + ic_{2})(c_{1} + ic_{2})/4 = c_{1}^2 - c_{2}^2 + i(c_{1}c_{2} + c_{2}c_{1}) = 0+i0 = 0.$$
Similarly one can
mathematically make two Majorana operators from any single Fermion operator via
$$c_{1} = (\psi + \psi^{\dagger})/2$$
$$c_{2} = (\psi - \psi^{\dagger})/(2i).$$
This simple relationship between the Fermion creation and annihilation algebra and an underlying Clifford algebra has long been a subject of speculation in physics. Only recently have experiments shown (indirect) 
evidence  \cite{Kouwenhouven} for Majorana Fermions underlying the electron.\\

\noindent {\bf Braiding.} Let there be given a set of Majorana operators
$\{ c_1, c_2, c_3, \ldots , c_n \}$ so that  $c_{i}^2 = 1$ for all $i$ and $c_{i}c_{j} = - c_{j}c_{i}$ for $i \ne j.$
Then there are natural braiding operators \cite{Ivanov,Kitaev} that act on the vector space with
these $c_k$ as the basis. The operators are mediated by algebra elements
$$\tau_{k} =(1 + c_{k+1} c_{k})/\sqrt{2},$$
$$\tau_{k}^{-1} = (1 - c_{k+1} c_{k})/\sqrt{2}.$$
Then the braiding operators are
$$T_{k}: Span \{c_1,c_2,\ldots, ,c_n \} \longrightarrow Span \{c_1,c_2,\ldots, ,c_n \}$$
via 
$$T_{k}(x) = \tau_{k} x \tau_{k}^{-1}.$$
The braiding is simply:
$$T_{k}(c_{k}) = c_{k+1},$$
$$T_{k}(c_{k+1}) = - c_{k},$$
and $T_{k}$ is the identity otherwise.
This gives a very nice unitary representaton of the Artin braid group and
it deserves better understanding.
\bigbreak

\noindent That there is much more to this braiding is indicated by the following result.\\

\noindent {\bf Clifford Braiding Theorem.} Let $C$ be the Clifford algebra over the real numbers generated by linearly independent elements $\{ c_{1},c_{2}, \ldots c_{n}\}$ with 
$c_{k}^2 = 1$ for all $k$ and $c_{k}c_{l} = - c_{l}c_{k}$ for $k \ne l.$
Then the  algebra elements $\tau_{k} =(1 + c_{k+1} c_{k})/\sqrt{2},$ form a representation of the (circular) Artin braid group.
That is, we have $\{\tau_{1},\tau_{2}, \ldots \tau_{n-1}, \tau_{n} \}$ where $\tau_{k} =(1 + c_{k+1} c_{k})/\sqrt{2}$ for $1 \le k < n$ and $\tau_{n} =(1 + c_{1} c_{n})/\sqrt{2},$
and $\tau_{k}\tau_{k+1}\tau_{k} = \tau_{k+1}\tau_{k}\tau_{k+1}$ for all $k$ and  $\tau_{i}\tau_{j} = \tau_{j}\tau_{i}$ when $|i-j|>2.$  Note that each braiding generator
$\tau_{k}$ has order $8.$ \\

\noindent {\bf Proof.} Let $a_{k} = c_{k+1}c_{k}.$
Examine the following calculation:
$$\tau_{k}\tau_{k+1}\tau_{k} = (\sqrt{2}/2)(1 + a_{k+1})(1 + a_{k})(1 + a_{k+1})$$
$$ = (\sqrt{2}/2)(1 + a_{k}  + a_{k+1} + a_{k+1}a_{k})(1 + a_{k+1})$$
$$ = (\sqrt{2}/2)(1 + a_{k}  + a_{k+1} + a_{k+1}a_{k} + a_{k+1} + a_{k}a_{k+1}  + a_{k+1}a_{k+1} + a_{k+1}a_{k}a_{k+1})$$
$$ = (\sqrt{2}/2)(1 + a_{k}  + a_{k+1} + c_{k+2}c_{k} + a_{k+1} + c_{k}c_{k+2}  -1 -c_{k}c_{k+1})$$
$$ = (\sqrt{2}/2)(a_{k}  + a_{k+1} + a_{k+1} + c_{k+1}c_{k})$$
$$ = (\sqrt{2}/2)(2a_{k}  + 2a_{k+1})$$
$$ = (\sqrt{2}/2)(a_{k}  + a_{k+1}).$$
Since the end result is symmetric under the interchange of $k$ and $k+1,$ we conclude that 
$$\tau_{k}\tau_{k+1}\tau_{k} = \tau_{k+1}\tau_{k}\tau_{k+1}.$$
Note that this braiding relation works circularly if we define $\tau_{n} =(1 + c_{1} c_{n})/\sqrt{2}.$
It is easy to see that $\tau_{i}\tau_{j} = \tau_{j}\tau_{i}$ when $|i-j|>2.$ This completes the proof. $\Box$\\

This representation of the (circular) Artin braid group is significant for the topological physics of Majorana Fermions.
This part of the structure needs further study. We discuss its relationship with the work of Mo-Lin Ge in the next section.\\

\noindent {\bf Remark.} It is worth noting that a triple of Majorana Fermions say $x,y,z$ gives rise
to a representation of the quaternion group. This is a generalization of
the well-known association of Pauli matrices and quaternions.
We have $x^2 = y^2 = z^2 = 1$ and, when different, they anti-commute.
Let $I = yx, J = zy, K = xz.$
Then $$I^2 = J^2 = K^2 = IJK = -1,$$ giving the quaternions.
The operators
$$X= \sigma_{I} = (1/\sqrt{2})(1 + I)$$
$$Y= \sigma_{J} = (1/\sqrt{2})(1 + J)$$
$$Z= \sigma_{K} = (1/\sqrt{2})(1 + K)$$
braid one another: $$XYX = YXY, XZX= ZXZ, YZY = ZYZ.$$
This is a special case of the braid group representation described above
for an arbitrary list of Majorana Fermions.
These braiding operators are entangling and so can be used for universal
quantum computation, but they give only partial topological quantum
computation due to the interaction with single qubit operators not
generated by them.
\bigbreak

\noindent {\bf Remark.} We can see just how the braid group relation arises in the Clifford Braiding Theorem as follows. Let $A$ and $B$ be given algebra elements with 
$$A^2 = B^2 = -1$$
and
$$AB = - BA.$$ 
(For example we can let $A = c_{k+1}c_{k}$ and B = $c_{k+2}c_{k+1}$ as above.)\\

Then we see immediately that 
$ABA = -BAA = -B(-1) = B$
and 
$BAB = A.$\\

\noindent Let $$\sigma_{A} = (1 + A)/\sqrt{2}$$ and $$\sigma_{B} = (1 + B)/\sqrt{2}.$$ 
Then $$\sigma_{A}\sigma_{B}\sigma_{A} = (1 +A)(1+B)(1 +A)/(2\sqrt{2})$$
$$ = (1+ B + A + AB)(1+A)/(2\sqrt{2})$$
$$= (1+ B + A + AB + A + BA + A^2 + ABA)/(2\sqrt{2})$$
$$ = (B + 2A + ABA)/(2\sqrt{2})$$
$$ = 2(A +B)/(2\sqrt{2})$$
$$ = (A +B)/(\sqrt{2}).$$
From this it follows by symmetry that 
$\sigma_{A}\sigma_{B}\sigma_{A} = \sigma_{B}\sigma_{A}\sigma_{B}.$\\

The relations $ABA = B$ and $BAB = A$ play a key role in conjunction with the fact that $A$ and $B$ have square equal to minus one.
We further remark that there is a remarkable similarity in the formal structure of this derivation and the way elements of the Temperley-Lieb algebra \cite{AnyonicTop,KnotLogic} can be combined to form a representation of the Artin Braid group.
In Temperley-Lieb algebra one has elements $R$ and $S$ with $RSR = R$ and $SRS = S$ where $R^{2} = R$ and $S^2 = S.$ The Temperley-Lieb algebra also gives rise to representations of the braid group and this is
used in the Fibonacci Model \cite{KnotLogic}. Both the Clifford algebra and the Temperley-Lieb algebra are associated with Majorana Fermions. The Clifford algebra is associated with the annihilation/creation algebra and the Temperley-Lieb algebra is associated with the fusion algebra (see \cite{KnotLogic} for a detailed discussion of the two algebras in relation to Majoranas). We will explore the relationship further in the next section.\\

Note that if we have a collection of operators $M_{k}$ such that 
$$M_{k}^{2} = -1,$$ 
$$M_{k}M_{k+1} = - M_{k+1}M_{k}$$ and
$$M_{i}M_{j} = M_{j}M_{i}$$ when $|i-j| \ge 2,$ then the operators
$$\sigma_{i} = (1 + M_{i})/\sqrt{2}$$ give a representation of the Artin Braid group. This is a way to generalize the Clifford Braiding Theorem, since defining $M_{i} = c_{i+1}c_{i}$ gives exactly such a set of operators, and 
the observation about braiding we give above extends to a proof of the more general theorem. One says that the operators $M_{k}$ generate an {\it extraspecial $2$-group} \cite{Rowell}.\\

\noindent {\bf Extraspecial 2-Group Braiding Theorem.} Let $\{ M_{1}, \cdots, M_{n} \}$ generate an extraspecial $2$-group as above. Then the operators $\sigma_{i} = (1 + M_{i})/\sqrt{2}$ for $i = 1,\cdots , n-1$ give a representation of the Artin Braid Group  with $\sigma_{i}\sigma_{j} = \sigma_{j}\sigma_{i}$ for $|i-j| \ge 2$ and $\sigma_{i}\sigma_{i+1}\sigma_{i} = \sigma_{i+1}\sigma_{i}\sigma_{i+1}.$ \\

\noindent {\bf Proof.} The proof follows at once from the discussion above. $\Box$\\

\noindent{\bf Remark.} In the next section we will see examples of braid group representations that arise from the Braiding Theorem, but are of a different character than the representations that come from the Clifford Braiding Theorem.\\

\begin{figure}
     \begin{center}
     \begin{tabular}{c}
     \includegraphics[height=5cm]{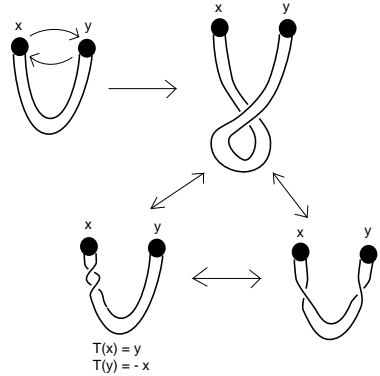}
     \end{tabular}
     \end{center}
     \caption{\bf Braiding Action on a Pair of Fermions}
     \label{ex}
     \end{figure}

The braiding operators in the Clifford Braiding Theorem can be seen to act on the vector space over the complex numbers that is spanned by the Majorana Fermions $\{ c_1, c_2, c_3, \ldots , c_n \}.$
To see how this works, let $x = c_{k}$ and $y = c_{k+1}$ from the basis above. Let
$$s =  \frac{1 + yx}{\sqrt{2}},$$
$$T(p) = sps^{-1} = (\frac{1 + yx}{\sqrt{2}})p(\frac{1 - yx}{\sqrt{2}}),$$ 
and verify that 
$T(x) = y$ and $T(y) = -x.$ Now view Figure~\ref{ex} where we have illustrated a 
topological interpretation for the braiding of two Fermions. In the topological interpretation the two Fermions are connected by a flexible belt. On interchange, the belt becomes twisted  by $2 \pi.$ 
In the topological interpretation a twist of $2 \pi$ corresponds to a phase change of $-1.$ 
(For more information on this topological interpretation of $2 \pi$ rotation for Fermions, see \cite{KP}.)
Without a further choice it is not evident which particle of the pair should receive the phase change. The topology alone tells us only the relative change of phase between the two particles. The Clifford algebra for Majorana Fermions makes a specific choice in the matter by using the linear ordering of the Majorana operators $\{ c_1 \ldots c_n \},$  and in this way fixes the representation of the braiding. In saying this, we are just referring to
the calculational form above where with $x$ prior to $y$ in the ordering, the braiding operator assigns the minus sign to one of them and not to the other.\\
\bigbreak

     \bigbreak

A remarkable feature of this braiding representation of Majorana Fermions is that it applies to give a representation of the $n$-strand braid group $B_{n}$ for any row of $n$ Majorana Fermions. It is not restricted to the 
quaternion algebra. Nevertheless, we shall now examine the braiding representations of the quaternions. These representations are very rich and can be used in situations (such as Fibonacci particles) involving particles that are their own anti-particles (analogous to the Majorana Fermions underlying electrons). Such particles can occur in collectivities of electrons as in the quantum Hall effect. In such situations it is theorized that one can examine the local interaction properties of the Majorana particles and then the braidings associated to triples of them (the quaternion cases) can come into play very strongly. In the case of electrons in nano-wires, one at the present time must make do with long range correlations between ends of the wires and forgoe such local interactions. Nevertheless, it is the purpose of this paper to juxtapose the full story about three strand braid group representations of the quaternions in the hope that this will lead to deeper understanding of the possibilities for even the electronic Majorana Fermions.\\

 \section {Remark on the Temperley-Lieb Algebra and Braiding from Majorana Operators}
 Let $\{ c_{1},\cdots, c_{n}\}$ be a set of Majorana Fermion operators so that $c_{i}^2 = 1$ and $c_{i}c_{j} = - c_{j}c_{i}$ for $i \ne j.$ Recall \cite{KP} that the $n$-strand Termperley-Lieb algebra $TL_{n}$ is generated by 
 $\{ U_{i},\cdots, U_{n-1}\}$ with the relations ( for an appropriate scalar $\delta$)  $U_{i}^2 = \delta U_{i}$, $U_{i}U_{i \pm 1}U_{i} = U_{i}$ and $U_{i}U_{j} = U_{j}U_{i}$ for $|i-j|\ge 2.$
Here we give a Majorana fermion representation of Temperley-Lieb algebra that is related to the Ivanov braid group representation that we have discussed in prvevious sections of the paper.\\

First define $A$ and $B$ as $A=c_{i}c_{i+1}$, $B=c_{i-1}c_{i}$ so that  $A^{2}=B^{2}=-1$ and $AB = - BA.$ Note the following relations:
Let  $U=(1+iA)$ and $V=(1+iB).$ Then $U^{2}=2U$ and $V^{2}=2V$ and $UVU=V$ and  $VUV=U.$\\

{\it Thus a Majorana fermion representation of the Temperley-Lieb algebra is given by:}
\begin{align}
&U_{k}=\frac{1}{\sqrt{2}}(1+ic_{k+1}c_{k}),\\
&U_{k}^{2}=\sqrt{2}U_{k},\\
&U_{k}U_{k\pm 1}U_{k} = U_{k},\\
&U_{k}U_{j}= U_{j}U_{k} \quad for |k-j|\ge 2.
\end{align}
Hence we have a representation of the Temperley-Lieb algebra with loop value $\sqrt{2}.$ Using this representation of the Temperley-Lieb algebra, we can construct (via the Jones representation of the braid group to the Temperley-Lieb agebra) another apparent represention of the braid group that is based on Majorana Fermions. In order to see this representation more explicitly, note that the Jones representation can be formulated as
$$\sigma_{k} = A U_{K} + A^{-1}1$$
$$\sigma_{k}^{-1} = A^{-1}U_{k} + A1$$
where $$-A^{2} - A^{-2} = \delta$$ and $\delta$ is the loop value for the Temperley-Lieb algebra so that $U_{k}^{2} = \delta U_{k}$ \cite{KP}.\\ 

In the case of the Temperley-Lieb algebra constructed from Majorana Fermion operators described above, we take $A$ such that $$-A^2 -A^{-2} = \sqrt{2}$$ and 
$$\sigma_{k} = A U_{k} + A^{-1}1$$
$$ = A \frac{1}{\sqrt{2}}(1+ic_{k+1}c_{k}) + A^{-1}1$$ 
$$=(A \frac{1}{\sqrt{2}}+ A^{-1}) +(i A \frac{1}{\sqrt{2}})c_{k+1}c_{k}.$$
Note that with $$A^2 + A^{-2} = - \sqrt{2},$$
$$x = (A \frac{1}{\sqrt{2}}+ A^{-1}), $$ and 
$$y = (i A \frac{1}{\sqrt{2}}),$$ then $x^2 = y^2.$
In fact the condition $x^2 = y^2$ is by itself a sufficient condition for the braiding relation to be satisfied.\\   

\noindent {\bf Remark.} Note that $x^2 - y^2 = (x + y)(x-y)$ and so the condition $x^2 = y^2$ is equivalent to the the condition that $x = + y$ or $x= -y.$ This means that the the braiding representatation that we have constructed for Majorana Fermion operators via the Temperley-Lieb representation  is not different from the one we first constructed from the Clifford Braiding Theorem. More precisely, we have the result of the Theorem below.\\

\noindent {\bf Theorem.} Suppose that $\alpha^2 = \beta^2 = -1$ and $\alpha \beta = - \beta \alpha.$ Define
$$\sigma_{\alpha} = x + y \alpha$$ and 
$$\sigma_{\beta} = x + y \beta.$$ Then $x^2 = y^2$ is a sufficient condition for the braiding relation
$$\sigma_{\alpha}\sigma_{\beta}\sigma_{\alpha} = \sigma_{\beta}\sigma_{\alpha}\sigma_{\beta} $$ to be satisfied.
Since $x^2 = y^2$ is satisfied by $y=x$ or $y = -x,$ we have the two possibilities:
$$\sigma_{\alpha} = x\sqrt{2}[(1 +  \alpha)/\sqrt{2}]$$ and 
$$\sigma_{\alpha} = x\sqrt{2}[(1 -  \alpha)/\sqrt{2}].$$
In other words, braiding representations that arise in this way will be equivalent to the original braiding represenatation of Ivanov that we have derived from the Clifford Braiding Theorem.\\

\noindent {\bf Proof.} Given the hypothesis of the Theorem, note that $\alpha \beta \alpha = -\beta \alpha^2 = \beta.$ Thus, we have
$$\sigma_{\alpha}\sigma_{\beta}\sigma_{\alpha} = (x + y \alpha)(x + y \beta)(x + y \alpha)$$
$$=(x^2 + xy\beta + xy \alpha + y^2 \alpha \beta)(x + y \alpha)$$
$$=x^3 + x^2 y\beta + x^2 y \alpha + xy^2 \alpha \beta + x^2y \alpha + xy^2 \beta \alpha + xy^2  \alpha^2 + y^3 \alpha \beta \alpha$$
$$=x^3 + x^2 y\beta + x^2 y \alpha  + x^2y \alpha  + xy^2  (-1) + y^3 \beta$$
$$=x^3 - xy^2  +( x^2 y + y^3) \beta + 2x^2 y \alpha.$$
This expression will be symmetric in $\alpha$ and $\beta$ if $x^2 y + y^3 = 2x^2 y,$ or equivalently, if $y(- x^2 + y^2) = 0.$ Thus the expression will be symmetric if
$x^2 = y^2.$ Since the braiding relation asserts the symmetry of this expression, this completes the proof of the Theorem. $\Box$\\

\section{\bf Majorana Fermions Generate Universal Braiding Gates}
Recall that in Section 3 we showed how to construct braid group representations. Let  $T_{k}: V_{n} \longrightarrow V_{n}$
defined by $$T_{k}(v) = \tau_{k} v \tau_{k}^{-1}$$ be defined as in Section 3.  Note that $ \tau_{k}^{-1} = \frac{1}{\sqrt{2}}(1 - c_{k+1} c_{k}).$ It is then easy to verify that
$$T_{k}(c_{k}) = c_{k+1},$$ $$T_{k}(c_{k+1}) = - c_{k}$$ and that $T_{k}$ is the identity otherwise.\\

For universality, take $n = 4$ and regard each $T_{k}$ as operating on $V \otimes V$ where $V$ is a single qubit space.  Then the braiding operator $T_{2}$ satisfies the Yang-Baxter equation and is an entangling operator. So we have universal gates (in the presence of single qubit unitary operators) from Majorana Fermions. If experimental work shows that Majorana Fermions can be detected and controlled, then it is possible that quantum computers based on these topological unitary representations will be constructed. Note that the matrix form $R$ of $T_{2}$ is
$$R = \left( \begin{array}{cccc}
1 & 0& 0 &0 \\
0 & 0& -1 &0 \\
0 & 1& 0 &0 \\
0 & 0& 0 &1\\
\end{array} \right)$$
Here we take the ordered basis $\{ |00\rangle, |01\rangle, |10\rangle, |11\rangle \}$ for the corresponding $2$-qubit space $V \otimes V$ so that 
$$R|00\rangle = |00\rangle  , R|01\rangle = |10\rangle,$$ $$R|10\rangle = -|01\rangle, R|11\rangle= |11\rangle.$$
It is not hard to verify that $R$ satisfies the Yang-Baxter Equation. To see that it is entangling we take the state $|\phi \rangle = a|0 \rangle + b|1\rangle$
and test $R$ on $$|\phi \rangle \otimes |\phi \rangle = a^2 |00\rangle +ab |01\rangle + ab |10\rangle + b^2 |11\rangle$$ and find that 
$$R(|\phi \rangle \otimes |\phi \rangle) = a^2 |00\rangle +ab |10\rangle  - ab |01\rangle + b^2 |11\rangle.$$
The determinant of this state is $a^2 b^2 + (ab)(ab) = 2a^2 b^2.$ Thus when both $a$ and $b$ are non-zero, we have that $R(|\phi \rangle \otimes |\phi \rangle)$ is entangled.
This proves that $R$ is an entangling operator, as we have claimed. This calculation shows that a fragment of the Majorana operator braiding can be used to make a universal quantum gate, and so to produce 
partial topological quantum computing if realized physically.\\

\section{\bf The Bell-Basis Matrix and Majorana Fermions}
We can say more about braiding  by using the operators $\tau_{k} = \frac{1}{\sqrt{2}}(1 + c_{k+1} c_{k}),$ as these operators have natural matrix representations. In particular, consider the 
Bell-Basis Matrix $B_{II}$  that is given as follows:
\begin{equation}
B_{II}=\frac{1}{\sqrt{2}}\left[\begin{array}{cccc}
1 & 0 & 0 & 1\\
0 & 1 & -1 & 0\\
0 & 1 & 1 & 0\\
-1 & 0 & 0 & 1
\end{array}\right]=\frac{1}{\sqrt{2}}\bigl(I+M\bigr)\quad\bigl(M^{2}=-1\bigr)
\end{equation}
where 
\begin{equation}
M=\left[\begin{array}{cccc}
0 & 0 & 0 & 1\\
0 & 0 & -1 & 0\\
0 & 1 & 0 & 0\\
-1 & 0 & 0 & 0
\end{array}\right]
\end{equation}
and we define $$M_{i} = I \otimes I \otimes \cdots I \otimes M \otimes I \otimes I \otimes \cdots \otimes I$$ where there are $n$ tensor factors in all and $M$ occupies the $i$ and $i+1$-st positions.
Then one can verify that theses matrices satisfy the relations of an ``extraspecial 2-group" \cite{Rowell,Ge}. The relations are as follows.
\begin{eqnarray}
M_{i}M_{i\pm1}&=&-M_{i\pm1}M_{i},\quad M^{2}=-I,\\
M_{i}M_{j}&=&M_{j}M_{i,}\quad \big|i-j\big|\geq2.
\end{eqnarray}

In fact, we can go further and identify a collection of Majorana Fermion operators behind these relations. Note that 

\begin{equation}
M=\left[\begin{array}{cccc}
0 & 0 & 0 & 1\\
0 & 0 & -1 & 0\\
0 & 1 & 0 & 0\\
-1 & 0 & 0 & 0
\end{array}\right] =
\left[\begin{array}{cccc}
1 & 0 & 0 & 0\\
0 & -1 & 0 & 0\\
0 & 0 & 1 & 0\\
0 & 0 & 0 & -1
\end{array}\right]
\left[\begin{array}{cccc}
0 & 0 & 0 & 1\\
0 & 0 & 1 & 0\\
0 & 1 & 0 & 0\\
1 & 0 & 0 & 0
\end{array}\right].
\end{equation}
Let
\begin{equation}
A =
\left[\begin{array}{cccc}
1 & 0 & 0 & 0\\
0 & -1 & 0 & 0\\
0 & 0 & 1 & 0\\
0 & 0 & 0 & -1
\end{array}\right], 
B= \left[\begin{array}{cccc}
0 & 0 & 0 & 1\\
0 & 0 & 1 & 0\\
0 & 1 & 0 & 0\\
1 & 0 & 0 & 0
\end{array}\right].
\end{equation}
Then we have $M= AB$ and $A^2 = B^2 = I$ while $AB + BA = 0.$ Thus we can take $A$ and $B$ themselves as Majorana Fermion operators.
By the same token, we have $$M_{i} =  I \otimes I \otimes \cdots I \otimes AB \otimes I \otimes I \otimes \cdots \otimes I$$ so that if we define
$$A_{i} =   I \otimes I \otimes \cdots I \otimes A \otimes I \otimes I \otimes \cdots \otimes I$$ and
$$B_{i} =   I \otimes I \otimes \cdots I \otimes B \otimes I \otimes I \otimes \cdots \otimes I.$$ 
Thus $$M_{i} = A_{i}B_{i}.$$

By direct calculation (details omitted here) , we find the following relations:
\begin{eqnarray}
A_{i}^2 =B_{i}^2 = 1,\\
A_{i}B_{i} = - B_{i}A_{i},\\
A_{i}B_{i+1} = - B_{i+1}A_{i},\\
A_{i+1}B_{i} = B_{i}A_{i+1},\\
A_{i}B_{j} = B_{j}A_{i}, \,\, for \,\, |i-j|>1 \\
A_{i}A_{j} = A_{j}A_{i} \\
B_{i}B_{j} = B_{j}B_{i} \,\, for \, all \,\, i \,\, and \,\, j.
\end{eqnarray}

\noindent {\bf Definition.} Call a list of pairs of Majorana operators that satisfies the above identities a {\it string} of Majorana operators or {\it Majorana string}.\\

\noindent {\bf Theorem.} Let $A_{i}, B_{i}$ ($i = 1, \cdots , n$) be a Majorana string. Then the product operators $M_{i} = A_{i}B_{i}$ ($i = 1, \cdots , n$) define an
extraspecial $2$-group and hence the operators $\sigma_{i} = (1 + M_{i})/\sqrt{2}$ for $i = 1,\cdots, n-1$ give a representation of the Artin Braid Group.\\

\noindent{\bf Proof.} We must verify that the $M_{i} = A_{i}B_{i}$ satisfy the relations for an extraspecial  $2$-group. To this end, note that 
$M_{i}^2 = A_{i}B_{i}A_{i}B_{i} = - A_{i}A_{i}B_{i}B_{i} = -1$ and that  (using the relations for a string of Majorana operators as above) $$M_{i}M_{i+1} = A_{i}B_{i}A_{i+1}B_{i+1} = A_{i}A_{i+1}B_{i}B_{i+1}$$ while
$$M_{i+1}M_{i} = A_{i+1}B_{i+1}A_{i}B_{i} = -  A_{i+1}A_{i}B_{i+1}B_{i}.$$ Thus $M_{i}M_{i+1} = - M_{i+1}M_{i}.$ We leave to the reader to check that for $|i-j|\ge 2,$
$M_{i}M_{j} = M_{j}M_{i}.$ This completes the proof that the $M_{i}$ form an extraspecial $2$-group. The braiding representation then follows from the Extraspecial 2-Group Braiding Theorem.$\Box$\\

\noindent {\bf Remark.} It is remarkable that these matrices representing Majorana Fermion operators give, so simply, a braid group representation and that the pattern behind this representation generalizes to any string of Majorana operators. Kauffman and Lomonaco \cite{BG} observed that $B_{II}$ satisfies the Yang-Baxter equation
and is an entangling gate. Hence $B_{II} = \frac{1}{\sqrt{2}}\bigl(I+M\bigr)\quad\bigl(M^{2}=-1\bigr)$ is a universal quantum gate in the sense of this paper.
It is of interest to understand the possible relationships of topological entanglement (linking and knotting) and quantum entanglement. See \cite{BG} for more than one point of view on this
question.\\

\noindent{\bf Remarks.}  The operators $M_{i}$ take the place here of the products of Majorana Fermions $c_{i+1}c_{i}$ in the Ivanov picture of braid group representation in the form
$\tau_{i} = (1/\sqrt{2})(1 + c_{i+1}c_{i}).$ This observation gives a concrete interpretation of these braiding operators and relates them to a Hamiltonian for a physical system by an observation of 
Mo-Lin Ge \cite{Ge}. Mo-Lin Ge shows that the observation of Ivanov \cite{Ivanov} that $\tau_{k} = (1/\sqrt{2})(1 + c_{k+1}c_{k}) = exp(c_{k+1}c_{k} \pi/4)$ can be extended
by defining $$\breve{R}_{k}(\theta) = e^{\theta c_{k+1} c_{k} }.$$ 
Then $\breve{R}_{i}(\theta)$ satisfies the 
full Yang-Baxter equation with rapidity parameter $\theta.$ That is, we have the equation 
$$\breve{R}_{i}(\theta_{1})\breve{R}_{i+1}(\theta_{2})\breve{R}_{i}(\theta_{3}) = \breve{R}_{i+1}(\theta_{3})\breve{R}_{i}(\theta_{2})\breve{R}_{i+1}(\theta_{1}).$$ 
This makes if very clear that $\breve{R}_{i}(\theta)$ has physical significance, and suggests examining the physical process for a temporal evolution of the unitary operator  $\breve{R}_{i}(\theta).$\\

In fact. following \cite{Ge},  we can construct a Kitaev chain \cite{Kitaev,KitaevFault} based on the solution $\breve{R}_i(\theta)$ of the Yang-Baxter Equation. Let a unitary evolution be governed by $\breve{R}_k(\theta)$. 
When $\theta$ in the unitary operator 
$\breve{R}_k(\theta)$ is time-dependent, we define a  state $|\psi(t)\rangle$ by $|\psi(t)\rangle=\breve{R}_k|\psi(0)\rangle$. With the Schroedinger equation 
$i\hbar\frac{\partial}{\partial t}|\psi(t)\rangle=\hat{H}(t)|\psi(t)\rangle$ one obtains:
\begin{equation}
i \hbar \frac{\partial}{\partial t}[\breve{R}_k|\psi(0)\rangle]=\hat{H}(t)\breve{R}_k|\psi(0)\rangle.
\end{equation}
Then the Hamiltonian $\hat{H}_k(t)$ related to the unitary operator $\breve{R}_k(\theta)$ is obtained by the formula:
\begin{equation}\label{SchrodingerEquation}
\hat{H}_i(t)=\textrm{i}\hbar\frac{\partial\breve{R}_k}{\partial t}\breve{R}_{k}^{-1}.
\end{equation}
Substituting $\breve{R}_k(\theta)=\exp(\theta c_{k+1}c_{k})$ into equation (\ref{SchrodingerEquation}), we have:
\begin{equation}\label{2MFHamiltonian}
\hat{H}_k(t)=\textrm{i}\hbar\dot{\theta}c_{k+1}c_{k}.
\end{equation}
This Hamiltonian describes the interaction between $k$-th and $(k+1)$-th sites via the parameter $\dot{\theta}$. When $\theta=n \times \frac{\pi}{4}$, the unitary evolution corresponds to the braiding progress of two nearest Majorana Fermion sites in the system as we have described it above. Here $n$ is an integer and signifies the time of the braiding operation. We remark that it is interesting to examine this periodicity of the appearance of the 
topological phase in the time evolution of this Hamiltonian. (Compare with discussion in \cite{RK}.) For applications, one may consider processes that let the Hamiltonian take the the system right to one of these topological points and then this Hamiltonian cuts off. This goes beyond the work of Ivanov, who examines the representation on Majoranas obtained by conjugating by these operators. The Ivanov representation is of order two,
while this representation is of order eight.\\

Mo-Lin Ge points out that if we only consider the nearest-neighbour interactions between Majorana Fermions,  and extend equation (\ref{2MFHamiltonian}) to an inhomogeneous chain with $2N$ sites, the derived model is expressed as:
\begin{equation}\label{YBEKitaev}
\hat{H}=\textrm{i}\hbar\sum_{k=1}^{N}(\dot{\theta}_1c_{2k}c_{2k-1}+\dot{\theta}_2c_{2k+1}c_{2k}),
\end{equation}
with $\dot{\theta}_1$  and $\dot{\theta}_2$ describing odd-even and even-odd pairs, respectively.\\

The Hamiltonian derived from $\breve{R}_{i}(\theta(t))$ corresponding to the braiding of nearest Majorana Fermion sites is exactly the same as the $1D$ wire proposed by Kitaev \cite{Kitaev}, and $\dot{\theta}_1=\dot{\theta}_2$ corresponds to the phase transition point in the ``superconducting'' chain. By choosing different time-dependent parameters $\theta_1$ and $\theta_2$, one finds that the Hamiltonian $\hat{H}$ corresponds to different phases.
These observations of Mo-Lin Ge give physical substance and significance to the Majorana Fermion braiding operators discovered by Ivanov \cite{Ivanov}, putting them into a robust context of Hamiltonian evolution via the simple Yang-Baxterization  $\breve{R}_{i}(\theta) = e^{\theta c_{i+1} c_{i} }.$ \\

In \cite{BG}, Kauffman and Lomonaco observe that the Bell Basis Change Matrix $B_{II},$ is a solution to the Yang-Baxter equation. This solution can be seen as a $4 \times 4$ matrix representation for the operator $\breve{R}_{i}(\theta).$ One can ask whether there is relation between topological order, quantum entanglement and braiding. This is the case for the Kitaev chain where non-local Majorana modes are entangled and have
braiding structure.\\

As we have pointed out above, Ge makes the further observation, that the Bell-Basis Matrix $B_{II}$ can be used to construct an extraspecial 2-group and a braid representation that has the same properties as
the Ivanov braiding for the Kitaev chain. We have pointed out that this extraspecial 2-group comes from a {\it string} of Majorana matrix operators (described above). This suggests alternate physical interpretations that need to be investigated in a sequel to this paper.\\

\section {\bf Appendix on $SU(2)$ Representations of the Artin Braid Group}
This appendix is determines all the representations of the three strand Artin braid group $B_{3}$ to the special unitary group $SU(2)$ and
concomitantly to the unitary group $U(2).$ One regards the groups $SU(2)$ and $U(2)$ as acting on a single qubit, and so $U(2)$ is usually regarded as the
group of local unitary transformations in a quantum information setting. If one is looking for a coherent way to represent all unitary transformations by
way of braids, then $U(2)$ is the place to start. Here we will give an example of a representation of the three-strand braid group
that generates a dense subset of $SU(2)$ (the Fibonacci model \cite{Kitaev,AnyonicTop}). Thus it is a fact that local unitary transformations can be "generated by braids" in many ways.
The braid group representations related to Majorana Fermions are also seen to have their roots in these quaternion representations, as we shall see at the end of the appendix.
\bigbreak

We begin with the structure of $SU(2).$ A matrix in $SU(2)$ has the form 
$$ M = 
\left( \begin{array}{cc}
z & w \\
-\bar{w} & \bar{z} \\
\end{array} \right),$$ where $z$ and $w$ are complex numbers, and $\bar{z}$ denotes the complex conjugate of $z.$ 
To be in $SU(2)$ it is required that $Det(M)=1$ and that $M^{\dagger} = M^{-1}$ where $Det$ denotes determinant, and $M^{\dagger}$ is the conjugate transpose of $M.$
Thus if
$z = a + bi$ and $w = c + di$ where $a,b,c,d$ are real numbers, and $i^2 = -1,$ then 
$$ M = 
\left( \begin{array}{cc}
a + bi & c + di \\
-c + di & a - bi \\
\end{array} \right)$$  with $a^2 + b^2 + c^2 + d^2 = 1.$ It is convenient to write
$$M =
a\left( \begin{array}{cc}
1 & 0 \\
0 & 1 \\
\end{array} \right) +
b\left( \begin{array}{cc}
i & 0\\
0 & -i \\
\end{array} \right) +
c\left( \begin{array}{cc}
0  & 1 \\
-1  & 0\\
\end{array} \right) +
d\left( \begin{array}{cc}
0 & i \\
i & 0 \\
\end{array} \right),$$ and to abbreviate this decomposition as
$$M = a + bI +cJ + dK$$
where 
$$ 1 \equiv
\left( \begin{array}{cc}
1 & 0 \\
0 & 1 \\
\end{array} \right),
I \equiv
\left( \begin{array}{cc}
i & 0\\
0 & -i \\
\end{array} \right),
J \equiv
\left( \begin{array}{cc}
0  & 1 \\
-1  & 0\\
\end{array} \right),
K \equiv
\left( \begin{array}{cc}
0 & i \\
i & 0 \\
\end{array} \right)$$ so that 
$$I^2 = J^2 = K^2 = IJK = -1$$ and 
$$IJ = K, JK=I, KI = J$$
$$JI = -K, KJ = -I, IK = -J.$$
The algebra of $1,I,J,K$ is called the {\it quaternions} after William Rowan Hamilton who discovered this algebra prior to the discovery of 
matrix algebra. Thus the unit quaternions are identified with $SU(2)$ in this way. We shall use this identification, and some facts about 
the quaternions to find the $SU(2)$ representations of braiding. First we recall some facts about the quaternions. For a detailed exposition of 
basics for quaternions see \cite{KP}.\\

\begin{enumerate}
\item Note that if $q = a + bI +cJ + dK$ (as above), then $q^{\dagger} = a - bI - cJ - dK$ so that $qq^{\dagger} = a^2 + b^2 + c^2 + d^2 = 1.$
\item A general quaternion has the form $ q = a + bI + cJ + dK$ where the value of $qq^{\dagger} = a^2 + b^2 + c^2 + d^2,$ is not fixed to unity.
The {\it length} of $q$ is by definition $\sqrt{qq^{\dagger}}.$
\item A quaternion of the form $rI + sJ + tK$ for real numbers $r,s,t$ is said to be a {\it pure} quaternion. We identify the set of pure
quaternions with the vector space of triples $(r,s,t)$ of real numbers $R^{3}.$
\item Thus a general quaternion has the form $q = a + bu$ where $u$ is a pure quaternion of unit length and $a$ and $b$ are arbitrary real numbers.
A unit quaternion (element of $SU(2)$) has the additional property that $a^2 + b^2 = 1.$
\item If $u$ is a pure unit length quaternion, then $u^2 = -1.$ Note that the set of pure unit quaternions forms the two-dimensional sphere
$S^{2} = \{ (r,s,t) | r^2 + s^2 + t^2 = 1 \}$ in $R^{3}.$
\item If $u, v$ are pure quaternions, then $$uv = -u \cdot v + u \times v$$ whre $u \cdot v$ is the dot product of the vectors $u$ and $v,$ and 
$u \times v$ is the vector cross product of $u$ and $v.$ In fact, one can take the definition of quaternion multiplication as
$$(a + bu)(c + dv) = ac + bc(u) + ad(v) + bd(-u \cdot v + u \times v),$$ and all the above properties are consequences of this
definition. Note that quaternion multiplication is associative.
\item Let $g = a + bu$ be a unit length quaternion so that $u^2 = -1$ and $a = cos(\theta/2), b=sin(\theta/2)$ for a chosen angle $\theta.$
Define $\phi_{g}:R^{3} \longrightarrow R^{3}$ by the equation $\phi_{g}(P) = gPg^{\dagger},$ for $P$ any point in $R^{3},$ regarded as a pure quaternion.
Then $\phi_{g}$ is an orientation preserving rotation of $R^{3}$ (hence an element of the rotation group $SO(3)$). Specifically, $\phi_{g}$ is a rotation
about the  axis $u$ by the angle $\theta.$ The mapping $$\phi:SU(2) \longrightarrow SO(3)$$ is a two-to-one surjective map from the special unitary group to
the rotation group. In quaternionic form, this result was proved by Hamilton and by Rodrigues in the middle of the nineteeth century.
The specific formula for $\phi_{g}(P)$ is shown below:
$$\phi_{g}(P) = gPg^{-1} = (a^2 - b^2)P + 2ab (P \times u) + 2(P \cdot u)b^{2}u.$$
\end{enumerate}

We want a representation of the three-strand braid group in $SU(2).$ This means that we want a homomorphism $\rho: B_{3} \longrightarrow SU(2),$ and hence
we want elements $g = \rho(s_{1})$ and $h= \rho(s_{2})$ in $SU(2)$ representing the braid group generators $s_{1}$ and $s_{2}.$ Since $s_{1}s_{2}s_{1} =
s_{2}s_{1}s_{2}$ is the generating relation for $B_{3},$ the only requirement on $g$ and $h$ is that $ghg = hgh.$ We rewrite this relation as
$h^{-1}gh = ghg^{-1},$ and analyze its meaning in the unit quaternions.
\bigbreak

Suppose that $g = a + bu$ and $h=c + dv$ where $u$ and $v$ are unit pure quaternions so that $a^2 + b^2 = 1$ and $c^2 + d^2 = 1.$
then $ghg^{-1} = c +d\phi_{g}(v)$ and $h^{-1}gh = a + b\phi_{h^{-1}}(u).$ Thus it follows from the braiding relation that 
$a=c,$ $b= \pm d,$ and that $\phi_{g}(v) = \pm \phi_{h^{-1}}(u).$  However, in the case where there is a minus sign we have
$g = a + bu$ and $h = a - bv = a + b(-v).$ Thus we can now prove the following Theorem. 
\bigbreak

\noindent {\bf Theorem.} Let $u$ and $v$ be pure unit quaternions and $g = a + bu$ and $h=c + dv$ have unit length. Then (without loss of generality), the braid relation $ghg=hgh$ is true if
and only if
$h = a + bv,$ and $\phi_{g}(v) = \phi_{h^{-1}}(u).$ Furthermore, given that $g = a +bu$ and $h = a +bv,$ the condition $\phi_{g}(v) = \phi_{h^{-1}}(u)$
is satisfied if and only if $u \cdot v = \frac{a^2 - b^2}{2 b^2}$ when $u \ne v.$ If $u = v$ then $g = h$ and the braid relation is trivially
satisfied.
\bigbreak

\noindent {\bf Proof.} We have proved the first sentence of the Theorem in the discussion prior to its statement. Therefore assume that
$g = a +bu, h = a +bv,$ and $\phi_{g}(v) = \phi_{h^{-1}}(u).$ 
We have already stated the formula for $\phi_{g}(v)$ in the discussion about quaternions:
$$\phi_{g}(v) = gvg^{-1} = (a^2 - b^2)v + 2ab (v \times u) + 2(v \cdot u)b^{2}u.$$ By the same token, we have
$$\phi_{h^{-1}}(u) = h^{-1}uh = (a^2 - b^2)u + 2ab (u \times -v) + 2(u \cdot (-v))b^{2}(-v)$$
$$= (a^2 - b^2)u + 2ab (v \times u) + 2(v \cdot u)b^{2}(v).$$ Hence we require that
$$(a^2 - b^2)v + 2(v \cdot u)b^{2}u = (a^2 - b^2)u + 2(v \cdot u)b^{2}(v).$$ This equation is equivalent to
$$2(u \cdot v)b^{2} (u - v) = (a^2 - b^2)(u - v).$$
If $u \ne v,$ then this implies that $$u \cdot v = \frac{a^2 - b^2}{2 b^2}.$$
This completes the proof of the Theorem. $\Box$\\
\bigbreak

\noindent{\bf The Majorana Fermion Example.} Note the case of the theorem where
$$g = a +bu, h = a +bv.$$ Suppose that $u \cdot v = 0.$ Then the theorem tells us that we need 
$a^2 - b^2 = 0$ and since $a^2 +b^2 = 1,$ we conclude that $a = 1/\sqrt{2}$ and $b$ likewise.
For definiteness, then we have for the braiding generators (since $I$, $J$ and $K$ are mutually orthogonal) the three operators
$$A = \frac{1}{\sqrt{2}}(1 + I),$$
$$B =\frac{1}{\sqrt{2}}(1 + J),$$
$$C = \frac{1}{\sqrt{2}}(1 + K).$$
Each pair satisfies the braiding relation so that $ABA = BAB, BCB = CBC, ACA =CAC.$ We have already met this braiding triplet in our discussion of the construction of braiding operators from Majorana Fermions in  Section 3. This shows (again) how close Hamilton's quaternions are to topology and how braiding is fundamental to the structure of Fermionic physics.
\bigbreak

\noindent{\bf The Fibonacci  Example.} Let
$$g = e^{I\theta} = a + bI$$ where $a = cos(\theta)$ and $b = sin(\theta).$
Let $$h = a + b[(c^2 - s^2)I + 2csK]$$ where $c^2 + s^2 = 1$ and $c^2 - s^2 = \frac{a^2 - b^2}{2b^2}.$ Then we can rewrite $g$ and $h$ in matrix form
as the matrices $G$ and $H.$ Instead of writing the explicit form of $H,$ we write $H = FGF^{\dagger}$ where $F$ is an element of $SU(2)$ as shown below.

$$G =
\left( \begin{array}{cc}
e^{i\theta} & 0 \\
0 & e^{-i\theta} \\
\end{array} \right)$$

$$F =
\left( \begin{array}{cc}
ic & is \\
is & -ic \\
\end{array} \right)$$
This representation of braiding where one generator $G$ is a simple matrix of phases, while the other generator $H = FGF^{\dagger}$ is derived from $G$ by
conjugation by a unitary matrix, has the possibility for generalization to representations of braid groups (with more than three strands) to $SU(n)$ or
$U(n)$ for 
$n$ greater than $2.$ In fact we shall see just such representations \cite{AnyonicTop} by using a version of topological quantum field theory.
The simplest example is given by 
$$g = e^{7 \pi I/10}$$
$$f = I\tau  + K \sqrt{\tau}$$
$$h = f g f^{-1}$$
where $\tau^{2} + \tau = 1.$
Then $g$ and $h$ satisfy $ghg=hgh$ and generate a representation of the three-strand braid group that is dense in $SU(2).$ We shall call this the 
{\it Fibonacci} representation of $B_{3}$ to $SU(2).$
\bigbreak

 At this point we can close this paper with the speculation that braid group representations such as this Fibonacci representation can be realized in the context of electrons in nano-wires.
 The formalism is the same as our basic Majorana representation. It has the form of a braiding operator of the form $$exp(\theta yx)$$ where $x$ and $y$ are Majorana operators and the angle 
 $\theta$ is not equal to $\pi/4$ as is required in the full Majorana representation. For a triple $\{x,y,z\}$ of Majorana operators, any quaternion representation is available. Note how this will effect the 
 conjugation representation: Let $T = r + s yx$ where $r$ and $s$ are real numbers with $r^2 + s^2 = 1$ (the cosine and sine of $\theta$), chosen so that a representation of the braid group is formed
 at the triplet (quaternion level). Then $T^{-1} = r - s yx$ and the reader can verify that 
 $$TxT^{-1} = (r^2 - s^2)x + 2rs y,$$
 $$TyT^{-1} =   (r^2 - s^2)y - 2rs x.$$
 Thus we see that the original Fermion exchange occurs with $r=s$ and then the sign on $-2rs$ is the well-known sign change in the exchange of Fermions. Here it is generalized to a more complex linear combination
 of the two particle/operators.\\


\begin{thebibliography}{99}


\bibitem{BA}  R.J. Baxter.  Exactly Solved Models in Statistical Mechanics.  Acad. Press (1982).

\bibitem{Beenakker} C. W. J. Beenakker, Search for Majorana Fermions in superconductors,
arXiv: 1112.1950.

\bibitem{B1}
N. E. Bonesteel, L. Hormozi, G. Zikos and S. H. Simon, Braid topologies for quantum computation,
 Phys. Rev. Lett. 95 (2005), no. 14, 140503, 4 pp. quant-ph/0505665.

\bibitem{B2}
S. H. Simon, N. E. Bonesteel, M. H. Freedman, N. Petrovic and L. Hormozi, Topological quantum computing with only one mobile quasiparticle,  Phys. Rev. Lett. 96 (2006), no. 7, 070503, 4 pp.,
quant-ph/0509175.


\bibitem{BB}
J. L. Brylinski and R. Brylinski, Universal quantum gates, in {\em Mathematics of Quantum Computation}, Chapman \& Hall/CRC Press, Boca Raton, Florida, 2002 (edited by R. Brylinski and G. Chen).




\bibitem{D}  
P.A.M. Dirac, {\em Principles of Quantum Mechanics},  Oxford University Press,
1958.


\bibitem{Fradkin}
E. Fradkin and P. Fendley, Realizing non-abelian statistics in time-reversal invariant systems, Theory Seminar, Physics Department, UIUC, 4/25/2005.

\bibitem{Ge} Li-Wei Yu and Mo-Lin Ge, More about the doubling degeneracy
operators associated with Majorana fermions and Yang-Baxter equation, Sci. Rep.\textbf{5},8102(2015).

\bibitem{Hatfield} B. Hatfield, ``Quantum Field Theory of Point Particles and Strings", Perseus Books, Cambridge, Massachusetts (1991).

\bibitem{Ivanov} 
D. A. Ivanov, Non-abelian statistics of half-quantum vortices in $p$-wave superconductors, 
Phys. Rev. Lett. 86, 268 (2001). 


\bibitem{KL}
L.H. Kauffman, {\em Temperley-Lieb Recoupling Theory and Invariants of Three-Manifolds},
Princeton University Press, Annals Studies {\bf 114} (1994). 


\bibitem {KP}
L.H. Kauffman, {\em Knots and Physics}, World Scientific Publishers (1991), 
Second Edition (1993), Third Edition (2002), Fourth Edition (2012).

\bibitem{KLogic} L. H. Kauffman, Knot logic and topological quantum computing with Majorana Fermions, in ``Logic and Algebraic Structures in Quantum Computing", edited by Chubb, Eskadarian and Harizanov,
Cambridge University Press (2016), pp. 223-335.

\bibitem{TEQE} 
L.H. Kauffman and S. J. Lomonaco Jr.,  Quantum entanglement and topological 
entanglement,  New Journal of Physics {\bf 4} (2002), 73.1--73.18 
(http://www.njp.org/).

\bibitem{Teleport}
L. H. Kauffman, Teleportation Topology, quant-ph/0407224, (in the Proceedings
of the  2004 Byelorus Conference on Quantum Optics),  {\it Opt. Spectrosc.} 9, 2005, 227-232.

 
\bibitem{Spie} 
L.H. Kauffman and S. J. Lomonaco Jr.,
Entanglement Criteria - Quantum and Topological, in 
{\em Quantum Information and Computation - Spie Proceedings, 21-22 April, 2003, 
Orlando, FL}, Donkor, Pinch and Brandt (eds.), Volume 5105, pp. 51--58.

\bibitem{QK1} L. H. Kauffman and S. J. Lomonaco Jr., Quantum knots, in {\it Quantum Information
and Computation II, Proceedings of Spie, 12 -14 April 2004} (2004), ed. by Donkor Pirich and Brandt, pp. 268-284.

\bibitem{QK2} S. J. Lomonaco and L.  H. Kauffman, Quantum Knots and Mosaics, 
Journal of Quantum Information Processing, Vol. 7, Nos. 2-3, (2008), pp. 85 - 115. 
http://arxiv.org/abs/0805.0339

 \bibitem{QK3} S. J. Lomonaco and L. H. Kauffman, Quantum knots and lattices, or a blueprint for quantum systems that do rope tricks. Quantum information science and its contributions to mathematics, 209Ð276, Proc. Sympos. Appl. Math., 68, Amer. Math. Soc., Providence, RI, 2010.
 
 \bibitem{QK4} S. J. Lomonaco and L. H. Kauffman, Quantizing braids and other mathematical structures: the general quantization procedure. In Brandt, Donkor, Pirich, editors, {\em Quantum Information and Comnputation IX - Spie Proceedings, April 2011}, Vol. 8057, of Proceedings of Spie, pp. 805702-1 to 805702-14, SPIE 2011.

\bibitem{QK5} L. H. Kauffman and S. J. Lomonaco, Quantizing knots groups and graphs. In Brandt, Donkor, Pirich, editors, {\em Quantum Information and Comnputation IX - Spie Proceedings, April
2011}, Vol. 8057, of Proceedings of Spie, pp. 80570T-1 to 80570T-15, SPIE 2011.

\bibitem{BG}
L. H. Kauffman and S. J. Lomonaco, Braiding Operators are Universal Quantum
Gates, New Journal of Physics 6 (2004) 134, pp. 1-39.

\bibitem{AnyonicTop}
Kauffman, Louis H.; Lomonaco, Samuel J., Jr. $q$-deformed spin networks, knot polynomials and anyonic topological quantum computation. J. Knot Theory Ramifications 16 (2007), no. 3, 267--332. 
 
\bibitem{SpinTop}
L. H. Kauffman and S. J. Lomonaco Jr.,
Spin Networks and Quantum Computation, in ``Lie Theory and Its Applications in Physics VII"
eds. H. D. Doebner and V. K. Dobrev, Heron Press, Sofia (2008), pp. 225 - 239.
 

\bibitem{QCJP1}
L. H. Kauffman, Quantum computing and the Jones polynomial, math.QA/0105255, in {\em Quantum Computation and Information}, S. Lomonaco, 
Jr. (ed.), AMS CONM/305, 2002, pp.~101--137.

\bibitem{QCJP2}
L. H. Kauffman and S. J. Lomonaco Jr.,
Quantum entanglement and topological entanglement. New J. Phys. 4 (2002), 73.1Ð73.18.

\bibitem{Fibonacci}
L. H. Kauffman and S. J. Lomonaco Jr., The Fibonacci Model and the Temperley-Lieb Algebra.
{\ it International J. Modern Phys. B}, Vol. 22, No. 29 (2008), 5065-5080.


\bibitem{Iterants}  Louis H Kauffman, Iterants, Fermions and Majorana Operators. In "Unified Field Mechanics - Natural Science Beyond the Veil of Spacetime", edited by
R. Amoroso, L. H. Kauffman, P. Rowlands, World Scientific Pub. Co. (2015).  pp. 1-32.

\bibitem{IAlg} Louis H. Kauffman, Iterant Algebra, Entropy 2017, 19, 347; doi:10.3390/e19070347 (30 pages).

\bibitem{KnotLogic} Louis H. Kauffman . Knot logic and topological quantum computing with Majorana Fermions. In ``Logic and algebraic structures in quantum computing and information", Lecture Notes in Logic, J. Chubb, J. Chubb, Ali Eskandarian, and V. Harizanov, editors,  124 pages Cambridge University Press (2016).

\bibitem{BMF} Louis H. Kauffman, Braiding and Majorana Fermions, Journal of Knot Theory and Its Ramifications Vol. 26, No. 9 (2017) 1743001 (21 pages).

\bibitem{KN1} L. H. Kauffman and P. Noyes, 
Discrete physics and the Dirac equation, Physics Lett. A, No.
218 (1996), pp. 139-146.



\bibitem{Kouwenhouven} 
V. Mourik,K. Zuo, S. M. Frolov, S. R. Plissard, E.P.A.M. Bakkers, L.P. Kouwenhuven, Signatures of Majorana Fermions in hybrid superconductor-semiconductor devices,
arXiv: 1204.2792.

 
\bibitem{Kitaev}
A. Kitaev, Anyons in an exactly solved model and beyond,  Ann. Physics 321 (2006), no. 1, 2Ð111.
{\em arXiv.cond-mat/0506438 v1 17 June 2005}.

\bibitem{KitaevFault} A. Kitaev, Fault-tolerant quantum computation by anyons,
Annals Phys. 303 (2003) 2-30. quant-ph/9707021.

\bibitem{Majorana} E. Majorana, A symmetric theory of electrons and positrons, 
I Nuovo Cimento,{\bf 14} (1937), pp. 171-184.

\bibitem{N}
M. A. Nielsen and I. L. Chuang, ``	Quantum Computation and Quantum Information," Cambrige University Press, Cambridge (2000). 

\bibitem{Rowell} J. Franko, E. C. Rowell, and Z. Wang, Extraspecial 2-Groups and Images of Braid Group Representations, 
JKTR Vol. 15, No. 4 (2006) 413Ð427, World Scientific Publishing Company.

\bibitem{RK} Rukhsan Ul Haq and L. H Kauffman, Z/2Z topological order and Majorana doubling in Kitaev Chain, (to appear)  arXiv:1704.00252v1 [cond-mat.str-el].


\bibitem{Wilczek}
F. Wilczek, {\em Fractional Statistics and Anyon Superconductivity,} World Scientific Publishing Company (1990).

\end{thebibliography}
\end{document}